\documentclass[11pt]{article}
\usepackage{amscd}
\usepackage{amsfonts}
\usepackage{amsmath}
\usepackage{amssymb}
\usepackage{amsthm}
\usepackage{bbm}
\usepackage{CJK}
\usepackage{fancyhdr}
\usepackage{graphicx}
\usepackage{hyperref}
\usepackage{indentfirst}
\usepackage{latexsym}
\usepackage{mathrsfs}
\usepackage{xypic}

\usepackage[top=1in,bottom=1in,left=1.25in,right=1.25in]{geometry}
\newtheorem{theorem}{Theorem}[section]
\newtheorem{lemma}[theorem]{Lemma}
\newtheorem{definition}[theorem]{Definition}
\newtheorem{proposition}[theorem]{Proposition}
\newtheorem{example}[theorem]{Example}

\newtheorem{corollary}[theorem]{Corollary}
\newtheorem{remark}[theorem]{Remark}

\usepackage[top=1in,bottom=1in,left=1.25in,right=1.25in]{geometry}
\def\<{\langle}
\def\>{\rangle}
\def\a{\alpha}
\def\b{\beta}

\def\c{\cdot}

\def\g{\gamma}

\def\o{\otimes}

\date{}
\begin{document}
\renewcommand{\baselinestretch}{1.2}
\renewcommand{\arraystretch}{1.0}
\title{\bf Some results on $\delta$-Hom-Jordan Lie conformal superalgebras}
\author{{\bf Shuangjian Guo$^{1}$,  Shengxiang Wang$^{2}$\footnote
        { Corresponding author(Shengxiang Wang):~~wangsx-math@163.com} }\\
{\small 1. School of Mathematics and Statistics, Guizhou University of Finance and Economics} \\
{\small  Guiyang  550025, P. R. of China} \\
{\small 2.~ School of Mathematics and Finance, Chuzhou University}\\
 {\small   Chuzhou 239000,  P. R. of China}}
 \maketitle
 \maketitle
\begin{center}
\begin{minipage}{13.cm}

{\bf \begin{center} ABSTRACT \end{center}}
 In this paper,  we introduce  the representation theory of $\delta$-Hom-Jordan Lie conformal superalgebras 
 and  discuss the cases of adjoint representations. Furthermore,  we develop the cohomology theory of Hom-Lie conformal superalgebras and discuss
some applications to the study of deformations of regular $\delta$-Hom-Jordan Lie conformal superalgebras. Finally, we introduce derivations of multiplicative
$\delta$-Hom-Jordan Lie conformal superalgebras and study their properties.


 \end{minipage}
 \end{center}
 \normalsize\vskip1cm

\section{INTRODUCTION}
\def\theequation{0. \arabic{equation}}
\setcounter{equation} {0}

Lie conformal superalgebras are introduced by Kac in \cite{Kac98}, in which he gave an axiomatic description of the singular part of the operator product expansion of chiral fields in conformal field theory. It is an useful tool to study vertex superalgebras and has many applications in the theory of Lie superalgebras.
 Moreover, Lie conformal superalgebras  have close connections to Hamiltonian formalism in the theory of nonlinear evolution equations.    Zhao, Yuan and Chen  developed deformation of Lie conformal superalgebras and   introduced derivations of multiplicative Lie conformal
superalgebras and studied their properties in \cite{Zhao2017}.

As a generalization of Lie superalgebras and Jordan Lie algebras, the notion of $\delta$-Jordan Lie superalgebra
$(L, [\c,\c]_L, \delta)$ was introduced in \cite{OK97}, which is intimately related to both Jordan-super and antiassociative algebras. The case of $\delta=1$ yields the Lie superalgebra, and   the other case of $\delta=-1$
is called a Jordan-Lie superalgebra because it turns out to be a Jordan superalgebra. It is   convenient to
consider both cases of $\delta=\pm 1$, and called $\delta$-Jordan Lie superalgebras. However, the results on $\delta$-Jordan
Lie algebras and $\delta$-Jordan Lie superalgebras are few in the present period.

In \cite{Ammar2010}, Ammar and Makhlouf introduced the notion of Hom-Lie superalgebras,
they gave a classification of Hom-Lie admissible superalgebras and proved a  graded version of Hartwig-Larsson-Silvestrov Theorem.
Later, Ammar, Makhlouf and Saadaoui \cite{Ammar2013} studied the representation and the cohomology of Hom-Lie superalgebras,
 and calculated the derivations and the second cohomology group of $q$-deformed Witt superalgebra.
In \cite{Yuan17}, Yuan introduced  Hom-Lie conformal superalgebra and proved that a Hom-Lie conformal superalgebra is equivalent to a Hom-Gel'fand-Dorfman superbialgebra.

  Recently, the Hom-Lie conformal algebra was introduced
and studied in \cite{Yuan14}, where it was proved that a Hom-Lie conformal algebra is equivalent to a Hom-Gel'fand-Dorfman bialgebra.   
Zhao, Yuan and Chen  developed the cohomology theory of Hom-Lie conformal algebras and discuss
some applications to the study of deformations of regular Hom-Lie conformal
algebras. Also, they  introduced  derivations of multiplicative Hom-Lie conformal
algebras and study their properties in \cite{Zhao2016}, which is different from \cite{Zhao2017}.  Ma and Chen introduced the notions of  $\delta$-Hom-Jordan Lie superalgebras and  discussed the concepts of $\a^{k}$-derivations, representations and $T^{\ast}$-extensions of $\delta$-Hom-Jordan Lie superalgebras in detail, and  established some cohomological characterizations in \cite{MC18}.

The following questions arise naturally:
1. How do we introduce the notion of a $\delta$-Hom-Jordan Lie conformal superalgebra?
2. How do we give  a representation theory of $\delta$-Hom-Jordan Lie conformal superalgebras?
3. How do we give  derivations of multiplicative $\delta$-Hom-Jordan Lie conformal superalgebras?

The aim of this article is to answer these questions.

Let us briefly describe the setup of the present article.
 In Section 2,  we introduce  a representation theory of $\delta$-Hom-Jordan Lie conformal superalgebras, in particular, we discuss the cases of adjoint representations. 
 In Section 3, we develop  cohomology theory of $\delta$-Hom-Jordan Lie conformal superalgebras and discuss
some applications to the study of deformations of regular $\delta$-Hom-Jordan Lie conformal superalgebras.
 In Section 4, we introduce the notion of derivations of multiplicative $\delta$-Hom-Jordan Lie conformal superalgebras and prove  the direct sum of the space of derivations is a $\delta$-Hom-Jordan Lie conformal superalgebra.

Throughout the paper, all algebraic systems are supposed to be over a field $\mathbb{C}$, and denote by $\mathbb{Z}_{+}$ the set of all nonnegative integers and by $\mathbb{Z}$.

\section{Representations of $\delta$-Hom-Jordan Lie conformal superalgebras }
\def\theequation{\arabic{section}. \arabic{equation}}
\setcounter{equation} {0}

In this section, we introduce  a representation theory of $\delta$-Hom-Jordan Lie conformal superalgebras and  discuss the cases of adjoint representations.

\begin{definition}
A $\delta$-Hom-Jordan Lie conformal superalgebra $R=R_{\bar{0}}\oplus R_{\bar{1}}$ is a  $\mathbb{Z}_2$-graded $\mathbb{C}[\partial]$-module equipped with an even linear endomorphism
$\a$ such that $\a\partial=\partial\a$, and  a $\mathbb{C}$-linear map
\begin{eqnarray*}
R\o R\rightarrow \mathbb{C}[\lambda]\o R,  ~~~a\o b\mapsto [a_\lambda b]
\end{eqnarray*}
such that $[R_{\phi\lambda}R_{\varphi}]\subseteq R_{\phi+\varphi}[\lambda]$, $\phi,\varphi\in \mathbb{Z}_2$, and the following axioms hold for $a,b,c\in R$
\begin{eqnarray}
&&[\partial a_{\lambda}b] =-\lambda[a_\lambda b],[ a_{\lambda}\partial b] =(\partial+\lambda)[a_\lambda b],\\
&&[a_{\lambda}b]=-\delta(-1)^{|a|||b|}[b_{-\lambda-\partial}a],\delta=\pm 1,\\
&&[\a(a)_\lambda[b_\mu c]]=\delta[[a_{\lambda}b]_{\lambda+\mu}\a(c)]+\delta(-1)^{|a|||b|}[\a(b)_{\mu}[a_{\lambda}c]].
\end{eqnarray}
\end{definition}
A $\delta$-Hom-Jordan Lie conformal superalgebra $(R, \delta, \a)$ is called finite if $R$ is a finitely generated $\mathbb{C}[\partial]$-module.
The rank of $R$ is its rank as a $\mathbb{C}[\partial]$-module.

A $\delta$-Hom-Jordan Lie conformal superalgebra $(R, \delta, \a)$  is called multiplicative if $\a$ is an algebra endomorphism, i.e., $\a([a_\lambda b])=[\a(a)_{\lambda}\a(b)]$ for any $a,b\in R$. In particular, if $\a$ is an algebra isomorphism, then $(R, \delta, \a)$  is regular.

\begin{remark}
(1) If $\delta=1$ in Definition 2.1, then the  $\delta$-Hom-Jordan Lie conformal superalgebra $(R, \delta, \a)$
is just the Hom-Lie conformal superalgebra in Yuan \cite{Yuan17}.

(2) If $\delta=1$ and $\a=id$ in Definition 2.1, then the  $\delta$-Hom-Jordan Lie conformal superalgebra $(R, \delta, \a)$
is just the Lie conformal superalgebra in Zhao \cite{Zhao2017}.
\end{remark}

\begin{example}
Let $(R,\delta, \a)$ be a $\delta$-Hom-Jordan Lie conformal superalgebra, and $\b: R\rightarrow R$ be an even endomorphism. Then $(R, \delta, \b\circ \a)$  ia a $\delta$-Hom-Jordan Lie conformal superalgebra, where $[x_{\lambda}y]=\b([x_{\lambda}y]), \forall x,y \in R$.

\end{example}
\begin{example}
Let $g=g_{\bar{0}}\oplus g_{\bar{1}}$ be a complex $\delta$-Jordan Lie superalgebra with Lie bracket $[-,-]$. Let $(Curg)_{\theta}:=\mathbb{C}[\partial]\o g_\theta$ be the free $\mathbb{C}[\partial]$-module. Then $Curg=(Curg)_{\bar{0}}\oplus (Curg)_{\bar{1}}$ is a $\delta$-Hom-Jordan Lie conformal superalgebra with $\lambda$-bracket given by
\begin{eqnarray*}
&& \a(f(\partial)\o a)=f(\partial)\o \a(a),\\
&&[(f(\partial)\o a)_{\lambda}(g(\partial)\o b)]=f(-\lambda)g(\partial+\lambda)\o [a,b], \forall a,b\in g.
\end{eqnarray*}
\end{example}

\begin{example}
Let $(R,\delta, \a)$ and  $(R',\delta, \a')$ be two $\delta$-Hom-Jordan Lie conformal superalgebras,   Then $(R\oplus R',\delta, \a+\b)$ is a $\delta$-Hom-Jordan Lie conformal superalgebra with
\begin{eqnarray*}
[(u_1+v_1)_\lambda (u_2+v_2)]=[u_{1 \lambda}v_1]+[u_{2 \lambda}v_2],\forall u_1,u_2\in R, v_1,v_2\in R',\\
(\a+\a')(u+v)=\a(u)+\a'(v), \forall u\in R, \forall v\in R'.
\end{eqnarray*}
\end{example}

\begin{definition}
Let $(M, \b)$ and $(N, \omega)$ be $\mathbb{Z}_2$-graded $\mathbb{C}[\partial]$-modules. A Hom-conformal linear map of degree $\theta$ from $M$ to $N$ is a sequence $f={f_{(n)}}_{n\in \mathbb{Z}_{\geq 0}}$ of $f_{(n)}\in Hom_{\mathbb{C}}(M, N)$ satisfying that
\begin{eqnarray*}
&&\partial_N f_{(n)}-f_{(n)}\partial_M=-n f_{(n-1)}, f_{\lambda}(M_\mu)\subseteq N_{\mu+\theta}, n\in \mathbb{Z}_{\geq 0}, \mu, \theta\in \mathbb{Z}_2,\\
&& \partial_M\b=\b\partial_M,  \partial_N\omega=\omega\partial_N, f_{(n)}\b=\omega f_{(n)}.
\end{eqnarray*}
Set $f_\lambda=\sum_{n=0}^{\infty}\frac{\lambda^{n}}{n!}f_{(n)}$. Then $f$ is a Hom-conformal linear map of degree $\theta$ if and only if
\begin{eqnarray*}
&&f_{\lambda} \partial_M=(\partial_N+\lambda)f_\lambda, f_{\lambda}(M_\mu)\subseteq N_{\mu+\theta}[\lambda],\\
&& \partial_M\b=\b\partial_M,  \partial_N\omega=\omega\partial_N, f_{\lambda}\b=\omega f_{\lambda}.
\end{eqnarray*}
\end{definition}

\begin{definition}
A $\delta$-Hom-associative conformal superalgebra $R$ is a  $\mathbb{Z}_2$-graded $\mathbb{C}[\partial]$-module equipped with an even linear endomorphism
$\a$ such that $\a\partial=\partial\a$, endowed with a $\lambda$-product from $R\o R$ to $\mathbb{C}[\partial]\o R$, for any $a,b,c\in R$, satisfying the following conditions:
\begin{eqnarray}
&&(\partial a)_\lambda b=-\lambda a_{\lambda}b, a_{\lambda}(\partial b)=(\partial+\lambda)a_{\lambda}b,\nonumber\\
&& \a(a)_{\lambda}(b_{\mu}c)=\delta(a_\lambda b)_{\lambda+\mu}\a(c).
\end{eqnarray}
\end{definition}

\begin{theorem}
Let $(R,\delta, \alpha)$ be a $\delta$-Hom-associative conformal superalgebra with an even linear endomorphism $\alpha$.
One can define
 \begin{eqnarray*}
[a_{\lambda}b]=a_{\lambda} b-\delta(-1)^{|a||b|}b_{-\lambda-\partial} a, \forall a,b\in R.
  \end{eqnarray*}
 Then $(R, \delta, \alpha)$ is a  $\delta$-Hom-Jordan Lie conformal superalgebra .
\end{theorem}

{\bf Proof.} We only prove (2.3). For this, we take $a,b,c\in R$ and calculate
\begin{eqnarray*}
&&[\a(a)_{\lambda}[b_{\mu}c]]\\
&=&[\a(a)_{\lambda}(b_{\mu}c-\delta(-1)^{|b||c|}c_{-\mu-\partial} b)]\\
&=&\a(a)_{\lambda}(b_{\mu}c)-\delta(-1)^{|a|(|b|+|c|)}(b_{\mu}c)_{-\lambda-\partial}\a(a)\\
&&-\delta(-1)^{|b||c|}\a(a)_{\lambda}(c_{-\mu-\partial} b)+(-1)^{|a||c|+|a||b|+|b||c|}(c_{-\mu-\partial} b)_{-\lambda-\partial}\a(a).
  \end{eqnarray*}
Similarly, we have
\begin{eqnarray*}
&&[[a_\lambda b]_{\lambda+\mu}\a(c)]\\
&=& [(a_{\lambda} b-\delta(-1)^{|a||b|}b_{-\lambda-\partial} a)_{_{\lambda+\mu}}\a(c)]\\
&=& (a_{\lambda} b)_{\lambda+\mu}\a(c)-\delta(-1)^{|c|(|a|+|b|)}\a(c)_{-\lambda-\mu-\partial}(a_{\lambda} b)\\
&&-\delta(-1)^{|a||b|}(b_{-\lambda-\partial} a)_{_{\lambda+\mu}}\a(c)+(-1)^{|a||c|+|a||b|+|b||c|}\a(c)_{-\lambda-\mu-\partial}(b_{-\lambda-\partial} a)_{\lambda+\mu}\\
&&(-1)^{|a||b|}[\a(b)_{\mu},[a_{\lambda}c]]\\
&=&(-1)^{|a||b|}\a(b)_{\mu}(a_{\lambda}c)-\delta(-1)^{|a|(|c|+|b|)}\a(b)_{\mu}(c_{-\lambda-\partial} a)\\
&&-\delta(-1)^{|a||c|}(a_{\lambda}c)_{-\mu-\partial}\a(b)+(-1)^{|a||b|}(c_{-\lambda-\partial} a)_{-\mu-\partial}\a(b).
  \end{eqnarray*}
By the associativity  (2.4), it is not hard to check that
\begin{eqnarray*}
[\a(a)_\lambda[b_\mu c]]=\delta[[a_{\lambda}b]_{\lambda+\mu}\a(c)]+\delta(-1)^{|a|||b|}[\a(b)_{\mu}[a_{\lambda}c]],
\end{eqnarray*}
as desired. And this finishes the proof.  \hfill $\square$

\begin{example}
Let $R=\mathbb{C}[\partial]e_1\oplus \mathbb{C}[\partial](e_2+e_3)$ be a free $\mathbb{Z}_2$-graded $\mathbb{C}[\partial]$-module and
 \[
 e_1=  \left(
   \begin{array}{ccc}
     0 & 0 & 1 \\
     0 & 0 & 0 \\
     0 & 0 & 0 \\
   \end{array}
 \right),
  e_2= \left (
   \begin{array}{ccc}
     0 & 1 & 0 \\
     0 & 0 & 0 \\
     0 & 0 & 0 \\
   \end{array}
 \right),
 e_3= \left(
   \begin{array}{ccc}
     0 & 0 & 0 \\
     0 & 0 & 1 \\
     0 & 0 & 0 \\
   \end{array}
 \right).
 \]
    Define\begin{eqnarray*}
       && \a(e_1)=\delta e_1, \a(e_2)=e_3, \a(e_3)=e_2,\\
&& [e_{1\lambda} e_1]= [e_{2\lambda} e_2]=[e_{3\lambda} e_3]=0, [e_{1\lambda} e_2]=[e_{1\lambda} e_3]=0, [e_{2\lambda}e_3]=\delta e_1.
          \end{eqnarray*}
One may check directly that $(R, \delta, \a)$ is a $\delta$-Hom-Jordan Lie conformal superalgebra.
\end{example}

Let Chom$(M, N)_{\theta}$ denote the set of Hom-conformal linear maps of degree $\theta$ from $M$ to $N$. Then  Chom$(M, N)$= Chom$(M, N)_{\bar{0}}\oplus$  Chom$(M, N)_{\bar{1}}$ is a $\mathbb{Z}_2$-graded $\mathbb{C}[\partial]$-module via
\begin{eqnarray*}
\partial f_{(n)}=-nf_{(n-1)},  ~~~\mbox{equivalently}~~~  \partial f_{\lambda}=-\lambda f_{\lambda}.
\end{eqnarray*}
The composition $f_\lambda g:L\rightarrow N\o \mathbb{C}[\lambda]$ of Hom-conformal linear maps $f: M\rightarrow N$ and $g: L\rightarrow M$ is given by
\begin{eqnarray*}
(f_\lambda g)_{\lambda+\mu}= f_\lambda g_\mu, \forall f,g\in  Chom(M, N).
\end{eqnarray*}
If $(M, \b)$ is a finitely generated $\mathbb{Z}_2$-graded $\mathbb{C}[\partial]$-module, then Cend$(M):=$Chom$(M, M)$ is an associative
conformal algebra with respect to the above composition. Thus, Cend$(M)$ becomes a
Hom-Lie conformal superalgebra, denoted as $gc(M)$, with respect to the following $\lambda$-bracket
\begin{eqnarray*}
[f_\lambda g] =f_\lambda g_{\mu-\lambda}-(-1)^{|f||g|}g_{\mu-\lambda}f_{\lambda}.
\end{eqnarray*}
\begin{definition}
Let $(R,\delta, \a)$ be a $\delta$-Hom-Jordan Lie conformal superalgebra. A representation of  $R$ is a triple $(\rho, M, \b)$, where $(M, \b)$ is a  $\mathbb{Z}_2$-graded $\mathbb{C}[\partial]$-module,  and   $\rho: R\rightarrow Cend{M}$ is an even  linear map  satisfying
\begin{eqnarray}
&&\rho(\partial(a))_{\lambda}=-\lambda\rho(a)_{\lambda},\rho\partial=\partial\rho, [\rho(a)_{\lambda}, \rho(b)_{\mu}]=\rho([a_\lambda b])_{\lambda+\mu},\nonumber\\
&&\rho([a_{\lambda}b])_{\lambda+\mu}\b=\rho(\a(a))_{\lambda}\rho(b)_{\mu}-\delta(-1)^{|a||b|}\rho(\a(b))_{\mu}\rho(a)_{\lambda}.
\end{eqnarray}

\end{definition}

\begin{proposition}
Let $(R,\delta, \a)$ be a $\delta$-Hom-Jordan Lie conformal superalgebra  and $(\rho, M, \b)$  the representation of  $R$. 
Then the  direct sum $R\oplus M$ is a $\delta$-Hom-Jordan Lie conformal superalgebra, where the $\lambda$-bracket $[\c_{\lambda}\c]$ on $R\oplus M$ is 
\begin{eqnarray*}
[(a+u)_{\lambda}(b+v)]_{M}=[a_{\lambda}b]+\delta\rho(a)_{\lambda}v-(-1)^{|a||b|}\rho(b)_{-\partial-\lambda}u,
\end{eqnarray*}
 the twist map $\a+\b: R\oplus M\rightarrow R\oplus M$ is
\begin{eqnarray*}
(\a+\b)(a+u)=\a(a)+\b(u), 
\end{eqnarray*}
where $a,b\in R, u,v\in M.$
\end{proposition}
{\bf Proof.} Note that $R\oplus M$ is equipped with a $\mathbb{C}[\partial]$-module structure via
\begin{eqnarray*}
\partial(a+u)=\partial(a)+\partial(u), \forall a\in R, u\in M.
\end{eqnarray*}
It is easy to check that $(\a+\b)\circ \partial=\partial\circ (\a+\b)$. Now, we check that (2.1) holds, for any $a,b\in R$ and $u,v\in M$, we have
\begin{eqnarray*}
[\partial(a+u)_{\lambda}(b+v)]_{M}&=&[(\partial a+\partial u)_{\lambda}(b+v)]_{M}\\
&=& [\partial a_{\lambda}b]+\rho(\partial a)_{\lambda}v-(-1)^{|a||b|}\rho(b)_{-\partial-\lambda}\partial u\\
&=& -\lambda[a_{\lambda}b]-\lambda\rho(a)_{\lambda}v-(-1)^{|a||b|}(\partial-\lambda-\partial)\rho(b)_{-\partial-\lambda}u\\
&=&  -\lambda([a_{\lambda}b]+\rho(a)_{\lambda}v-(-1)^{|a||b|}\rho(b)_{-\partial-\lambda}u)\\
&=&  -\lambda [(a+u)_{\lambda}(b+v)]_{M},
\end{eqnarray*}
and
\begin{eqnarray*}
[(a+u)_{\lambda}\partial(b+v)]_{M}&=& [(a+u)_{\lambda}(\partial b+\partial v)]_{M}\\
&=& [a_\lambda \partial b]+\rho( a)_{\lambda}\partial v-(-1)^{|a||b|}\rho(\partial b)_{-\partial-\lambda}u\\
&=&(\partial+\lambda) [a_\lambda  b]+(\partial+\lambda)\rho( a)_{\lambda} v-(-1)^{|a||b|}(\partial+\lambda)\rho(b)_{-\partial-\lambda}u\\
&=& (\partial+\lambda)([a_{\lambda}b]+\rho(a)_{\lambda}v-(-1)^{|a||b|}\rho(b)_{-\partial-\lambda}u)\\
&=& (\partial+\lambda)[(a+u)_{\lambda}(b+v)]_{M}.
\end{eqnarray*}
For (2.2), we have
\begin{eqnarray*}
[(b+v)_{-\partial-\lambda}(a+u)]_{M}&=&[b_{-\partial-\lambda}a]+\rho(b)_{-\partial-\lambda}u-(-1)^{|a||b|}\rho(a)_{\lambda}v\\
&=&-\delta(-1)^{|a||b|}([a_{\lambda}b]-(-1)^{|a||b|}\rho(b)_{-\partial-\lambda}u+\rho(a)_{\lambda}v)\\
&=& -\delta(-1)^{|a||b|}[(a+u)_{\lambda}(b+v)]_{M}.
\end{eqnarray*}
Next we verify that the Hom-Jacobi identity, we compute
\begin{eqnarray*}
&&[(\a+\b)(a+u)_{\lambda}[(b+v)_{\mu}(c+w)]_M]_M\nonumber\\
&=& [(\a(a)+\b(u))_{\lambda}[(b+v)_{\mu}(c+w)]_M]_M\nonumber\\
&=& [(\a(a)+\b(u))_{\lambda}[[b_{\mu}c]+\rho(b)_{\mu}w-(-1)^{|b||c|}\rho(c)_{-\partial-\mu}v]_M\nonumber\\
&=& [\a(a)_{\lambda}[b_{\mu}c]]+\rho(\a(a))_{\lambda}(\rho(b)_{\mu}w)-\delta(-1)^{|c||b|}\rho(\a(a))_{\lambda}(\rho(c)_{-\partial-\mu}v)\nonumber\\
&&-(-1)^{(|c|+|b|)|a|}\rho([b_\mu c])_{-\partial-\lambda}(\b(u)),\\
&&(-1)^{|a||b|}[(\a+\b)(b+v)_{\mu}[(a+u)_{\lambda}(c+w)]_M]_M\nonumber\\
&=& (-1)^{|a||b|}[\a(b)_{\mu}[a_{\lambda}c]]+(-1)^{|a||b|}\rho(\a(b))_{\mu}(\rho(a)_{\lambda}w)\nonumber\\
&&-\delta(-1)^{|a||b|+|c||a|}\rho(\a(b))_{\mu}(\rho(c)_{-\partial-\lambda}u)-(-1)^{|c||b|}\rho([a_{\lambda}c]_{-\partial-\mu})\b(v),\\
&&[[(a+u)_{\lambda}(b+v)]_{M_{\lambda+\mu}}(\a+\b)(c+w)]_M\nonumber\\
&=& [[a_{\lambda}b]_{\lambda+\mu}\a(c)]+ \rho([a_{\lambda}b])_{\lambda+\mu}\b(w)-\delta(-1)^{|c|(|a|+|b|)}\rho(\a(c))_{-\partial-\lambda-\mu}(\rho(a)_{\lambda}v)\nonumber\\ &&+(-1)^{|c|(|b|+|a|)+|b||a|}\rho(\a(c))_{-\partial-\lambda-\mu}(\rho(b)_{-\partial-\lambda}u)~~~~~~~
\end{eqnarray*}

Since $(\rho, M, \b)$ is a representation of the $\delta$-Hom-Jordan Lie conformal superalgebra $R$, we have
\begin{eqnarray*}
\rho(\a(a))_{\lambda}(\rho(b)_{\mu}w)-\delta(-1)^{|a||b|}\rho(\a(b))_{\mu}(\rho(a)_{\lambda}w)=\rho([a_{\lambda}b])_{\lambda+\mu}\b(w).
\end{eqnarray*}
Thus  $(R\oplus M, \a+\b)$ is a $\delta$-Hom-Jordan Lie conformal superalgebra. \hfill $\square$

\begin{proposition}
Let $(R,\delta, \a)$ be a $\delta$-Hom-Jordan Lie conformal superalgebra. Define a Hom-conformal linear map $ad: R\rightarrow  Cend{R}$   by $(ad (a))_{\lambda}(b)=\delta[a_{\lambda}b]$.
Then $(R, \delta, ad, \a)$ is a representation of $R$.
\end{proposition}
{\bf Proof.} Since $(R,\delta, \a)$ is a $\delta$-Hom-Jordan Lie conformal superalgebra, the Hom-Jacobi identity (2.3) may be written as follows
\begin{eqnarray*}
ad([a_\lambda b])_{\lambda+\mu}(\a(c))=ad(\a(a))_{\lambda}(ad(b)_{\mu}(c))-\delta(-1)^{|a||b|}ad(\a(b))_{\mu}(ad(a)_{\lambda}(c)),
\end{eqnarray*}
 for any $a,b,c\in R$.
Then the conformal linear map $ad$ satisfies
\begin{eqnarray*}
ad([a_\lambda b])_{\lambda+\mu} \a=ad(\a(a))_{\lambda} ad(b)_{\mu}-\delta(-1)^{|a||b|}ad(\a(b))_{\mu} ad(a)_{\lambda}.
\end{eqnarray*}
Therefore, $(R, \delta,  ad, \a)$ is a representation of $R$.   \hfill $\square$

We call the representation $ad$ defined above the adjoint representation of the $\delta$-Hom-Jordan Lie conformal superalgebra.

In the following, we will explore the dual representation and coadjoint representation of  $\delta$-Hom-Jordan Lie conformal superalgebras.

Let $(R,\delta, \a)$ be a $\delta$-Hom-Jordan Lie conformal superalgebra and $(\rho, M, \b)$   a representation of $R$. Set $M^{\ast}$ be the dual vector space of $M$. We define a Hom-conformal linear map $\tilde{\rho}: R\rightarrow Cend{M^{\ast}}$ by $\tilde{\rho}_{\lambda}(a)=-\rho_{\lambda}(x)$.

Let $f\in M^{\ast}, a,b\in R$ and $u\in M$. We compute the right hand side of (2.4)
\begin{eqnarray*}
&&(\tilde{\rho}(\a(a))_{\lambda})\tilde{\rho}(b)_{\mu}-\delta(-1)^{|a||b|}(\tilde{\rho}(\a(b))_{\mu})\tilde{\rho}(a)_{\lambda} (f)(u)\\
&=& (\tilde{\rho}(\a(a))_{\lambda})(\tilde{\rho}(b)_{\mu}(f))-\delta(-1)^{|a||b|}(\tilde{\rho}(\a(b))_{\mu})(\tilde{\rho}(a)_{\lambda}(f)) (u)\\
&=&  -(\tilde{\rho}(b)_{\mu}(f))(\rho(\a(a))_{\lambda}(u))+\delta(-1)^{|a||b|}(\tilde{\rho}(a)_{\lambda}(f))(\tilde{\rho}(\a(b))_{\mu} (u))\\
&=& f(\rho(b)_{\mu}\rho(\a(a))_{\lambda}(u))-\delta(-1)^{|a||b|}f(\rho(a)_{\lambda}\rho(\a(b))_{\mu}(u))\\
&=& f((\rho(b)_{\mu}\rho(\a(a))_{\lambda}-\delta(-1)^{|a||b|}f(\rho(a)_{\lambda}\rho(\a(b))_{\mu})(u)).
\end{eqnarray*}
On the other hand, we set that the map $\tilde{\b}=\b$, then the left  hand side of (2.4),
\begin{eqnarray*}
((\tilde{\rho}([a_\lambda b])_{\lambda+\mu}\tilde{\b}) (f))(u)=((\tilde{\rho}([a_\lambda b])_{\lambda+\mu} (f\b))(u))=-f \b(\rho([a_\lambda b])_{\lambda+\mu}(u)).
\end{eqnarray*}
Therefore, we have the following proposition:
\begin{proposition}
Let $(R,\delta, \a)$ be a $\delta$-Hom-Jordan Lie conformal superalgebra and $(\rho, M, \b)$ the representation of $R$.
Then $(M^{\ast},\tilde{\rho}, \tilde{\b} )$ is a representation of    $(R,\delta, \a)$ if and only if
\begin{eqnarray*}
\b\rho([a_\lambda b])_{\lambda+\mu}=(-1)^{|a||b|}\rho(a)_{\lambda}\rho(\a(b)_{\mu}-   \delta\rho(b)_{\mu}\rho(\a(a))_{\lambda}, \forall a,b\in R.
\end{eqnarray*}
\end{proposition}

\begin{corollary}
Let $(R,\delta, \a)$ be a $\delta$-Hom-Jordan Lie conformal superalgebra and $(R, \delta, ad, \a)$   the adjoint representation of $R$, where $ad: g\rightarrow Cend{g}$. We set $\tilde{ad}: g\rightarrow Cend{g^{\ast}}$ and $\tilde{ad}(x)_{\lambda}(f)=-f ad(x)_{\lambda}$. Then $(g^{\ast},\tilde{ad}, \tilde{\a} )$  is a representation of $R$ if and only if
\begin{eqnarray*}
\a([a_{\lambda}b]_{\lambda+\mu}, c)=\delta(-1)^{|a||b|}[a_{\lambda} [\a(b)_{\mu}c]]-[b_{\mu}[\a(a)_\lambda c]], \forall a,b,c\in R.
\end{eqnarray*}
\end{corollary}

\section{ Nijenhuis operators of $\delta$-Hom-Jordan Lie conformal superalgebras }
\def\theequation{\arabic{section}. \arabic{equation}}
\setcounter{equation} {0}

In this section, we introduce the notions of Nijenhuis operators and deformations  of $\delta$-Hom-Jordan Lie conformal superalgebras  and show that the deformation generated
by a 2-cocycle Nijenhuis operator is trivial.

\begin{definition}
Let $(R,\delta, \a)$ be a $\delta$-Hom-Jordan Lie conformal superalgebra.
An $n$-cochain ($n\in \mathbb{Z}_{\geq0}$) of a regular $\delta$-Hom-Jordan Lie conformal superalgebra $R$  with coefficients in a module $(M,\b )$ is a $\mathbb{C}$-linear map
\begin{eqnarray*}
&&\g:R^n\rightarrow M[\lambda_1,... ,\lambda_n], \\
&&(a_1,...,a_n)\mapsto \g_{\lambda_1,...,\lambda_n}(a_1,... , a_n),
\end{eqnarray*}
where $M[\lambda_1,...,\lambda_n]$ denotes the space of polynomials with coefficients in $M$, satisfying the following conditions:

(1) Conformal antilinearity:
\begin{eqnarray*}
\g_{\lambda_1,...,\lambda_n}(a_1,... ,\partial a_i,... ,a_n)=-\lambda_i\g_{\lambda_1,...,\lambda_n}(a_1,... ,a_i,... ,a_n).
\end{eqnarray*}

(2) Skew-symmetry:
\begin{eqnarray*}
&&\gamma(a_{1},..., a_{i}, a_{j},..., a_{n})=-\delta(-1)^{|a_i||a_j|}\gamma(a_{1}, ..., a_{j}, a_{i}, ... , a_{n}).
\end{eqnarray*}

(3)Commutativity:
\begin{eqnarray*}
\g\circ \a=\beta\circ \g
\end{eqnarray*}
\end{definition}

Let $R^{\otimes 0} = \mathbb{C}$ as usual so that a $0$-cochain
 is an element of $M$. Define a differential $d^n: C^n_{\a}(R, M)\rightarrow C^{n+1}_{\a}(R, M)(n=1,2)$ by
 \begin{eqnarray*}
d^1 \g_{\lambda_1,\lambda_2}(u_1,u_2)&=&(-1)^{|\g||u_1|}\rho (\a(u_1))_{\lambda_1}\g(u_2)-\delta(-1)^{(|\g|+|u_1|)|u_2|}\rho(\a(u_2))_{\lambda_2}\g(u_1)]-\delta\g([u_{1\lambda_1}u_2]),\\
d^2 \g_{\lambda_1,\lambda_2, \lambda_3}(u_1,u_2, u_3)&=&=(-1)^{|\g||u_1|}\rho(\a(u_1))_{\lambda_1}\g(u_{2\lambda_2 }u_3)-\delta(-1)^{(|\g|+|u_1|)|u_2|}\rho(\a^{2}(u_2))_{\lambda_2}\g(u_{1\lambda_1}u_3)\\
&& +(-1)^{(|\g|+|u_1|+|u_2|)|u_3|}\rho(\a^{2}(u_3))_{\lambda_2}\g(u_{1\lambda_1}u_2)-\g([u_{1\lambda_1}u_2],\a(u_3))\\
&& +(-1)^{(|u_1|+|u_2|)|u_3|}(-1)^{|u_1||u_3|}\delta \g([u_{1\lambda_1}u_3],\a(u_2))\\
&& -(-1)^{|u_1||u_2|} (-1)^{(|u_1|+|u_2|)|u_3|}(-1)^{|u_2||u_3|}\g([u_{2\lambda_2}u_3],\a(u_1)).
\end{eqnarray*}
where $\g$ is extended linearly over the polynomials in
$\lambda_1,\lambda_2, \lambda_3$. In particular, if $\g$ is a 0-cochain, then $(d\g)_\lambda a=a_\lambda \g$.

\begin{proposition}
With the above notations, we have $d^2\circ d^1=0$.
\end{proposition}

 {\bf Proof.} For any $u_1,u_2,u_3\in R$, we have
 \begin{eqnarray*}
&& d^2 \circ d^1 \g_{\lambda_1, \lambda_2, \lambda_3}(u_1,u_2,u_3)\\
&=&(-1)^{|\g||u_1|}\rho(\a(u_1))_{\lambda_1}d^1\g_{\lambda_2, \lambda_3}(u_{2\lambda_2 }u_3)-\delta(-1)^{(|\g|+|u_1|)|u_2|}\rho(\a^{2}(u_2))_{\lambda_2}d^1\g_{\lambda_1, \lambda_3}(u_{1\lambda_1}u_3)\\
&& +(-1)^{(|\g|+|u_1|+|u_2|)|u_3|}\rho(\a^{2}(u_3))_{\lambda_3}d^1\g_{\lambda_1, \lambda_2}(u_{1\lambda_1}u_2)-\g_{\lambda_1, \lambda_2, \lambda_3}([u_{1\lambda_1}u_2]_{\lambda_1+\lambda_2}\a(u_3))\\
&& +(-1)^{(|u_1|+|u_2|)|u_3|}(-1)^{|u_1||u_3|}\delta d^1\g_{\lambda_1, \lambda_2, \lambda_3}([u_{1\lambda_1}u_3]_{\lambda_1+\lambda_3}\a(u_2))\\
&& -(-1)^{|u_1||u_2|} (-1)^{(|u_1|+|u_2|)|u_3|}(-1)^{|u_2||u_3|}d^1\g([u_{2\lambda_2}u_3]_{ \lambda_2+\lambda_3}\a(u_1))\\
&=&   (-1)^{|\g||u_1|}\rho(\a(u_1))_{\lambda_1}\\
&&((-1)^{|\g||u_2|}\rho(\a(u_2))_{\lambda_2}\g(u_3)-\delta (-1)^{(|\g|+|u_2|)|u_3|}\rho(\a(u_3))_{\lambda_3}\g(u_1)- \delta\g([u_{2\lambda_2}u_3]))\\
&& -\delta(-1)^{(|\g|+|u_1|)|u_2|}\rho(\a^{2}(u_2))_{\lambda_2}\\
&&((-1)^{|\g||u_1|}\rho(\a(u_1))_{\lambda_1}\g(u_3)-\delta (-1)^{(|\g|+|u_1|)|u_3|}\rho(\a(u_3))_{\lambda_3}\g(u_1)- \delta\g([u_{1\lambda_1}u_3]))
 \end{eqnarray*}
 \begin{eqnarray*}
&&+ (-1)^{(|\g|+|u_1|+|u_2|)|u_3|}\rho(\a^{2}(u_3))_{\lambda_3}\\
&& ((-1)^{|\g||u_1|}\rho(\a(u_1))_{\lambda_1}\g(u_2)-\delta (-1)^{(|\g|+|u_1|)|u_2|}\rho(\a(u_2))_{\lambda_2}\g(u_1)- \delta\g([u_{1\lambda_1}u_2]))\\
&& -((-1)^{|\g|(|u_1|+|u_2|)}\rho(\a([u_{1\lambda_1}u_2]))_{\lambda_1+\lambda_2}\g(\a(u_3))-\delta (-1)^{(|\g|+|u_1|+|u_2|)|u_3|}\rho(\a^2(u_3))_{\lambda_3}\g([u_{1\lambda_1}u_2])\\
&&- \delta\g([[u_{1\lambda_1}u_2]_{\lambda_1+\lambda_2}\a(u_3)]))\\
&& +\delta (-1)^{|u_2||u_3|}\\
&&((-1)^{|\g|(|u_1|+|u_3|)}\rho(\a([u_{1\lambda_1}u_3]))_{\lambda_3}\g(\a(u_2))-\delta (-1)^{(|\g|+|u_1|+|u_3|)|u_2|}\rho(\a^2(u_2))_{\lambda_2}\g([u_{1\lambda_1}u_3])\\
&&- \delta\g([[u_{1\lambda_1}u_3]_{\lambda_1+\lambda_3}\a(u_2)]))\\
&& -(-1)^{|u_1||u_2|+|u_1||u_3|}\\
&& ((-1)^{|\g|(|u_2|+|u_3|)}\rho(\a([u_{2\lambda_2}u_3]))_{\lambda_2+\lambda_3}\g(\a(u_1))-\delta (-1)^{(|\g|+|u_2|+|u_3|)|u_1|}\rho(\a^2(u_1))_{\lambda_1}\g([u_{2\lambda_2}u_3])\\
&&- \delta\g([[u_{2\lambda_2}u_3]_{\lambda_2+\lambda_3}\a(u_1)])).
 \end{eqnarray*}
Because $(M, \b)$ is an $R$-module, we have
\begin{eqnarray*}
&& \rho(\a([u_{2\lambda_2}u_3]))_{\lambda_2+\lambda_3}\g(\a(u_1))\\
&=& \rho(\a([u_{2\lambda_2}u_3]))_{\lambda_2+\lambda_3}\b \circ \g(u_1)\\
&& \rho(\a^2(u_2))_{\lambda_2}(\rho(\a(u_3))_{\lambda_3})\g(u_1)-\delta(-1)^{|u_2||u_3|}\rho(\a^2(u_3))_{\lambda_3}(\rho(\a(u_2))_{\lambda_2} \g(u_1).
\end{eqnarray*}
Similarly, we have
\begin{eqnarray*}
&& \rho(\a([u_{1\lambda_1}u_3]))_{\lambda_1+\lambda_3}\g(\a(u_2))\\
&=& \rho(\a^2(u_1))_{\lambda_1}(\rho(\a(u_3))_{\lambda_3})\g(u_2)-\delta(-1)^{|u_1||u_3|}\rho(\a^2(u_3))_{\lambda_3}(\rho(\a(u_1))_{\lambda_1} \g(u_2),\\
&& \rho(\a([u_{1\lambda_1}u_2]))_{\lambda_1+\lambda_2}\g(\a(u_3))\\
&=& \rho(\a^2(u_1))_{\lambda_1}(\rho(\a(u_2))_{\lambda_2})\g(u_3)-\delta(-1)^{|u_1||u_2|}\rho(\a^2(u_2))_{\lambda_2}(\rho(\a(u_1))_{\lambda_1} \g(u_3).
\end{eqnarray*}
 It follows that   $d^2\circ d^1=0$, as desired. \hfill $\square$

 Let $R$ be a regular $\delta$-Hom-Jordan Lie conformal superalgebra. Define
 \begin{eqnarray}
\rho(a)_\lambda b=\delta[\a^{s}(a)_\lambda b],~~~\mbox{for any $a,b\in R$}.
 \end{eqnarray}

 Let $\g\in C^n_{\a}(R,R_s)$.  Define an operator $d_s: C^{n}_{\a}(R,R_s)\rightarrow C^{n+1}_{\a}(R,R_s) (n=1,2)$ by
\begin{eqnarray*}
d_1 \g_{\lambda_1,\lambda_2}(u_1,u_2)&=&(-1)^{|\g||u_1|}\delta[\a^{1+s}(u_1)_{\lambda_1}\g(u_2)]\\
&&-(-1)^{(|\g|+|u_1|)|u_2|}[\a^{1+s}(u_2)_{\lambda_2}\g(u_1)]-\delta\g([u_{1\lambda_1}u_2]),\\
d_2 \g_{\lambda_1,\lambda_2, \lambda_3}(u_1,u_2, u_3)&=&(-1)^{|\g||u_1|}\delta[\a^{1+s}(u_1)_{\lambda_1}\g(u_{2\lambda_2 }u_3)]\\
&&-(-1)^{(|\g|+|u_1|)|u_2|}[\a^{2+s}(u_2)_{\lambda_2}\g(u_{1\lambda_1}u_3)]\\
&& +(-1)^{(|\g|+|u_1|+|u_2|)|u_3|}[\a^{2+s}(u_3)_{\lambda_2}\g(u_{1\lambda_1}u_2)]-\g([u_{1\lambda_1}u_2],\a(u_3))\\
&& +(-1)^{(|u_1|+|u_2|)|u_3|}(-1)^{|u_1||u_3|}\delta \g([u_{1\lambda_1}u_3],\a(u_2))\\
&& -(-1)^{|u_1||u_2|} (-1)^{(|u_1|+|u_2|)|u_3|}(-1)^{|u_2||u_3|}\g([u_{2\lambda_2}u_3],\a(u_1)).
\end{eqnarray*}

Taking $s=-1$, because $\psi\in C^{2}_{\a}(R,R_{-1})_{\overline{0}}$ is a bilinear operator commuting with $\a$,  we consider a $t$-parameterized family of bilinear operations on $R$,
\begin{eqnarray}
[a_\lambda b]_t=[a_\lambda b]+t\psi_{\lambda, -\partial-\lambda}(a,b),~~~\forall a,b\in R.
\end{eqnarray}
If $[\c_{\lambda}\c]$ endows $(R,[\c_{\lambda}\c], \delta, \a)$ with a $\delta$-Hom-Jordan Lie conformal superalgebra, we say that $\psi$  generates
a deformation of the $\delta$-Hom-Jordan Lie conformal superalgebra $R$. It is easy to see that $[\c_{\lambda}\c]$ satisfies (2.1) and (2.2).

 If it is true for (2.3), by expanding the Hom-Jacobi identity for $[\c_{\lambda}\c]$, we have
 \begin{eqnarray*}
 &&[\a(a)_\lambda[b_{\mu}c]]+t([\a(a)_\lambda(\psi_{\mu,-\partial-\mu}(b,c)]+\psi_{\lambda,-\partial-\lambda}(\a(a),[b_{\mu}c]))\\
 &&+t^2\psi_{\lambda,-\partial-\lambda}(\a(a),\psi_{\mu,-\partial-\mu}(b,c))\\
 &=&(-1)^{|a||b|} [\a(b)_\mu[a_{\lambda}c]]+t([\a(b)_\mu(\psi_{\lambda,-\partial-\lambda}(a,c))]+\psi_{\mu,-\partial-\mu}(\a(b),[a_{\lambda}c]))\\
 &&+(-1)^{|a||b|}t^2\psi_{\mu,-\partial-\mu}(\a(b),\psi_{\lambda,-\partial-\lambda}(a,c))+[[a_{\lambda}b]_{\lambda+\mu}\a(c)]\\
 &&+t([(\psi_{\lambda,-\partial-\lambda}(a,b))_{\lambda+\mu}\a(c)]+ \psi_{\lambda+\mu,-\partial-\lambda-\mu}([a_\lambda b],\a(c)))\\
 &&+t^2\psi_{\lambda+\mu,-\partial-\lambda-\mu}(\psi_{\lambda,-\partial-\lambda}(a,b), \a(c)).
 \end{eqnarray*}
This is equivalent to the following conditions
\begin{eqnarray}
 &&\psi_{\lambda,-\partial-\lambda}(\a(a),\psi_{\mu,-\partial-\mu}(b,c))\nonumber\\
 &=&(-1)^{|a||b|}\psi_{\mu,-\partial-\mu}(\a(b),\psi_{\lambda,-\partial-\lambda}(a,c))+ \psi_{\lambda+\mu,-\partial-\lambda-\mu}(\psi_{\lambda,-\partial-\lambda}(a,b),\a(c)),\\
&&[\a(a)_\lambda(\psi_{\mu,-\partial-\mu}(b,c)]+\psi_{\lambda,-\partial-\lambda}(\a(a),[b_{\mu}c])\nonumber\\
&=& (-1)^{|a||b|}[\a(b)_\mu(\psi_{\lambda,-\partial-\lambda}(a,c))]+(-1)^{|a||b|}\psi_{\mu,-\partial-\mu}(\a(b),[a_{\lambda}c])\nonumber\\
&&+[(\psi_{\lambda,-\partial-\lambda}(a,b))_{\lambda+\mu}\a(c)]+ \psi_{\lambda+\mu,-\partial-\lambda-\mu}([a_\lambda b],\a(c)).
\end{eqnarray}

By conformal antilinearity of $\psi$, we have
\begin{eqnarray}
[(\psi_{\lambda, -\partial-\lambda}(a,b))_{\lambda+\mu}c]=[\psi_{\lambda, \mu}(a,b)_{\lambda+\mu}c].
\end{eqnarray}

On the other hand, Let $\psi$ be a cocycle, i.e., $d_{-1}\psi=0$. In fact,
\begin{eqnarray}
0&=&(d^{-1}\psi)_{\lambda, \mu, \g}(a,  b, c)\nonumber\\
&=& \delta[\a(a)_{\lambda}\psi_{\mu, \g}(b, c)]-(-1)^{|a||b|}[\a(b)_{\mu}\psi_{\lambda, \g}(a, c)]+\delta(-1)^{(|a|+|b|)|c|}[\a(c)_{\g}\psi_{\lambda, \mu}(a, b)]\nonumber\\
&& -\psi_{\lambda+\mu, \g}([a_\lambda b], \a(c))+\delta(-1)^{|b||c|}\psi_{\lambda+\g, \mu}([a_\lambda c], \a(b))-(-1)^{|a|(|b|+|c|)}\psi_{\g+\mu, \lambda}([b_\mu c], \a(a))\nonumber\\
&=& \delta[\a(a)_{\lambda}\psi_{\mu, \g}(b, c)]-\delta(-1)^{|a||b|}[\a(b)_{\mu}\psi_{\lambda, \g}(a, c)]-\delta\psi_{\lambda, \mu}(a, b)_{-\partial-\g}\a(c)]\nonumber\\
&& +\delta\psi_{\lambda, \g+\mu}(\a(a), [b_\mu c])-\delta(-1)^{|b||a|}\psi_{\mu, \lambda+\g}(\a(b), [a_\lambda c])-\delta\psi_{\lambda+\mu, \g}([a_\lambda b], \a(c)).
\end{eqnarray}

Replacing $\g$ by $-\lambda-\mu-\partial$ in (3.6), we have
\begin{eqnarray*}
&&[\a(a)_\lambda(\psi_{\mu,-\partial-\mu}(b,c)]+\psi_{\lambda,-\partial-\lambda}(\a(a),[b_{\mu}c])\nonumber\\
&=& (-1)^{|a||b|}[\a(b)_\mu(\psi_{\lambda,-\partial-\lambda}(a,c))]+(-1)^{|a||b|}\psi_{\mu,-\partial-\mu}(\a(b),[a_{\lambda}c])\nonumber\\
&&+[(\psi_{\lambda,-\partial-\lambda}(a,b))_{\lambda+\mu}\a(c)]+ \psi_{\lambda+\mu,-\partial-\lambda-\mu}([a_\lambda b],\a(c)).
\end{eqnarray*}

When $\psi$ is a 2-cocycle satisfying (3.3), $(R,\delta, \a)$ forms  a $\delta$-Hom-Jordan Lie conformal superalgebra. In this case, $\psi$ generates a deformation of the $\delta$-Hom-Jordan Lie conformal superalgebra $R$.

A deformation is said to be trivial if there is a linear operator $f\in \widetilde{C}_{\a}^{1}(R,R_{-1})$ such that for $T_t=id+tf$, there holds
\begin{eqnarray}
T_t([a_\lambda b]_t)=[T_t(a)_\lambda T_t(b)],~~~\forall a,b\in R.
\end{eqnarray}
\begin{definition}
 A linear operator $f\in \widetilde{C}_{\a}^{1}(R,R_{-1})$ is a Hom-Nijienhuis operator
if  \begin{eqnarray}
    [f(a)_\lambda f(b)]=f([a_\lambda b]_{N}),~~~\forall a,b\in R,
    \end{eqnarray}
where the bracket $[\c,\c]_{N}$ is defined by
\begin{eqnarray}
[a_\lambda b]_{N}=[f(a)_\lambda b]+[a_\lambda f(b)]-f([a_\lambda b]),~~~\forall a,b\in R.
\end{eqnarray}
\end{definition}
\begin{theorem}
 Let $(R,\delta, \a)$ be a regular $\delta$-Hom-Jordan Lie conformal superalgebra, and $f\in C_{\a}^{1}(R,R_{-1})$
a Hom-Nijienhuis operator. Then a deformation of $(R,\delta, \a)$) can be obtained by putting
\begin{eqnarray}
\psi_{\lambda, -\partial-\lambda}(a,b)=[a_\lambda b]_{N},~~~\forall a, b\in R.
\end{eqnarray}
Furthermore, this deformation is trivial.
\end{theorem}
{\bf Proof.}  To show that $\psi$  generates a deformation, we need to verify (3.3). By (3.9)
and (3.10), we have
\begin{eqnarray*}
&&\psi_{\lambda,-\partial-\lambda}(\a(a),\psi_{\mu,-\partial-\mu}(b,c))-(-1)^{|a||b|}\psi_{\mu,-\partial-\mu}(\a(b),\psi_{-\partial-\lambda,\lambda}(c,a))\\
&&-\psi_{-\partial-\lambda-\mu, \lambda+\mu}(\a(c), \psi_{\lambda,-\partial-\lambda}(a,b))\\
&=&[f(\a(a))_{\lambda}[f(b)_{\mu}c]]+[f(\a(a))_{\lambda}[b_\mu f(c)]]-[f(\a(a))_{\lambda}f([b_{\mu}c])] \\
&&+[\a(a)_{\lambda}[f(b)_\mu f(c)]]-f([\a(a)_{\lambda}[f(b)_\mu c]])-f([\a(a)_{\lambda}[b_\mu f(c)]])+f([\a(a)_{\lambda}f([b_\mu c])])\\
&&-(-1)^{|a||b|}[f(\a(b))_{\mu}[f(c)_{-\partial-\lambda}a]]-(-1)^{|a||b|}[f(\a(b))_{\mu}[c_{-\partial-\lambda} f(a)]]\\
&&+(-1)^{|a||b|}[f(\a(b))_{\mu}f([c_{-\partial-\lambda}a])]-(-1)^{|a||b|}[\a(b)_{\mu}[f(c)_{-\partial-\lambda} f(a)]]\\
&&+(-1)^{|a||b|}f([\a(b)_{\mu}[f(c)_{-\partial-\lambda} a]])+(-1)^{|a||b|}f([\a(b)_{\mu}[c_{-\partial-\lambda} f(a)]])\\
&&-(-1)^{|a||b|}[f(\a(b))_{\mu}f([c_{-\partial-\lambda} a])])-[f(\a(c))_{-\partial-\lambda-\mu}[f(a)_{\lambda}b]]-[f(\a(c))_{-\partial-\lambda-\mu}[a_\lambda f(b)]]\\
&&+[f(\a(c))_{-\partial-\lambda-\mu}f([a_{\lambda}b])]-[\a(c)_{-\partial-\lambda-\mu}[f(a)_{\lambda}f(b)]]\\
&&+f([\a(c)_{-\partial-\lambda-\mu}[f(a)_{\lambda}b]])+f([\a(c)_{-\partial-\lambda-\mu}[a_{\lambda}f(b)]])-f([\a(c)_{-\partial-\lambda-\mu}f([a_{\lambda}b])]).
\end{eqnarray*}
Since $f$ is a Hom-Nijenhuis operator, we get
\begin{eqnarray*}
&&-[f(\a(a))_{\lambda}f([b_{\mu}c])]+f([\a(a)_{\lambda}f([b_\mu c])])\\
&=&-f([f(\a(a))_{\lambda}[b_{\mu}c]])+f^2([\a(a)_{\lambda}[b_\mu c]]),\\
&& -(-1)^{|a||b|}[f(\a(b))_{\mu}f([a_{\lambda}c])]+(-1)^{|a||b|}f([\a(b)_{\mu}f([a_{\lambda} c])])\\
&=&-(-1)^{|a||b|}f([f(\a(b))_{\mu}[a_{\lambda}c]])+(-1)^{|a||b|}f^2([\a(b)_{\mu}[a_{\lambda} c]]),
\end{eqnarray*}
\begin{eqnarray*}
&& -[f(\a(c))_{-\partial-\lambda-\mu}f([a_{\lambda}b])]+f([c_{-\partial-\lambda-\mu}f([a_{\lambda}b])])\\
&=&-f[f(\a(c))_{-\partial-\lambda-\mu}[a_{\lambda}b]]+f^2([c_{-\partial-\lambda-\mu}[a_{\lambda}b]])
\end{eqnarray*}
Note that
\begin{eqnarray*}
[\a(a)_\lambda[f(b)_\mu f(c)]]=[[a_{\lambda}f(b)]_{\lambda+\mu}f(\a(c))]+(-1)^{|a||b|}[f(\a(b))_{\mu}[a_{\lambda} f(c)]].
\end{eqnarray*}
It follows that
\begin{eqnarray*}
 &&\psi_{\lambda,-\partial-\lambda}(\a(a),\psi_{\mu,-\partial-\mu}(b,c))\nonumber\\
 &=&(-1)^{|a||b|}\psi_{\mu,-\partial-\mu}(\a(b),\psi_{\lambda,-\partial-\lambda}(a,c))+ \psi_{\lambda+\mu,-\partial-\lambda-\mu}(\psi_{\lambda,-\partial-\lambda}(a,b),\a(c)).
\end{eqnarray*}
In the similar way, one can check that  (3.4). This proves that $\psi$ generates a deformation of the regular $\delta$-Hom-Jordan Lie conformal superalgebra $(R, \delta, \a)$.

Let $T_t=id+tf$. By  (3.2) and (3.10), we have
\begin{eqnarray}
T_t([a_\lambda b]_t)&=&(id+tf)([a_\lambda b]+t\psi_{\lambda,-\partial-\lambda}(a,b))\nonumber\\
&=& (id+tf)([a_\lambda b]+t[a_\lambda b]_N)\nonumber\\
&=& [a_\lambda b]+t([a_\lambda b]_N+f([a_\lambda b]))+t^2f([a_\lambda b]_N).
\end{eqnarray}
On the other hand, we have
\begin{eqnarray}
[T_t(a)_\lambda T_t(b)]&=&[(a+tf(a))_\lambda(b+tf(b))]\nonumber\\
&=& [a_\lambda b]+t([f(a)_\lambda b]+[a_\lambda f(b)])+t^2[f(a)_\lambda f(b)].
\end{eqnarray}
Combining  (3.11) with (3.12) gives $T_t([a_\lambda b]_t)=[T_t(a)_\lambda T_t(b)]$. Therefore the deformation
is trivial.\hfill $\square$

\section{Derivations of multiplicative  Hom-Lie conformal superalgebras }
\def\theequation{\arabic{section}. \arabic{equation}}
\setcounter{equation} {0}
In this section, we study derivations of multiplicative $\delta$-Hom-Jordan Lie conformal superalgebras and prove  the direct sum of the space of derivations is also a $\delta$-Hom-Jordan Lie conformal superalgebra.
\begin{definition} Let $(R,\delta, \a)$ be a multiplicative $\delta$-Hom-Jordan Lie conformal superalgebra. Then a
Hom-conformal linear map $D_\lambda: R\rightarrow R$ is called an $\a^k$-derivation of $(R,\a)$ if
\begin{eqnarray}
D_\lambda \circ \a=\a\circ D_\lambda,\nonumber\\
D_\lambda([a_\mu b])=\delta^{k}[D_\lambda(a)_{\lambda+\mu}\a^{k}(b)]+\delta^{k}(-1)^{|a||D|}[\a^{k}(a)_{\mu}D_\lambda(b)],
\end{eqnarray}
for any $a,b\in R.$
\end{definition}
Denote by $Der_{\a^{k}}$ the set of $\a^{k}$-derivations of the multiplicative $\delta$-Hom-Jordan Lie conformal superalgebra $(R,\delta, \a)$. For any $a\in R$ satisfying $\a(a)=a$, define $D_{k}:R\rightarrow R$ by
\begin{eqnarray*}
D_{k}(a)_{\lambda}(b)=\delta[a_\lambda\a^{k+1}(b)],~\delta^k=1,~~\forall b\in R.
\end{eqnarray*}

Then $D_{k}(a)$ is an $\a^{k+1}$-derivation, which is called an inner $\a^{k+1}$-derivation. In fact
\begin{eqnarray*}
D_{k}(a)_{\lambda}(\partial b)&=&\delta[a_\lambda\a^{k+1}(\partial b)]\\
&=& \delta[a_\lambda\partial\a^{k+1}( b)]\\
&=& (\partial+\lambda)D_{k}(a)_{\lambda}(b),\\
D_{k}(a)_{\lambda}(\a(b))&=&\delta[a_\lambda\a^{k+2}(b)]\\
&=& \delta\a[a_\lambda\a^{k+1}( b)]\\
&=& \a\circ D_{k}(a)_{\lambda}(b),\\
D_{k}(a)_{\lambda}([b_{\mu}c])&=&\delta [a_\lambda \a^{k+1}([b_{\mu}c])\\
&=& \delta [\a(a)_\lambda [\a^{k+1}(b)_{\mu}\a^{k+1}(c)]\\
&=& \delta[a_\lambda\a^{k+1}(b)]_{\lambda+\mu}\a^{k+1}(c)]+\delta(-1)^{|a||b|}[\a^{k+1}(b)_{\mu}[\a(a)_{\lambda}\a^{k+1}(c)]]\\
&=&\delta[a_\lambda\a^{k+1}(b)]_{\lambda+\mu}\a^{k+1}(c)]+\delta(-1)^{|a||b|}[\a^{k+1}(b)_{\mu}[a_{\lambda}\a^{k+1}(c)]]\\
&=&\delta^{k+1}[D_{k}(a)_{\lambda}(b)_{\lambda+\mu}\a^{k+1}(c)]+\delta^{k+1}(-1)^{|a||b|}[\a^{k+1}(b)_{\mu}(D_{k}(a)_{\lambda}(c))].
\end{eqnarray*}
Denote by $Inn_{\a^k} (R)$ the set of inner $\a^k$-derivations.
\begin{eqnarray*}
Inn^{\a^k}(R)=\{\delta[u_\lambda\a^{k+1}(\cdot)]|\delta^k=1,\a(u)=u~~\forall u\in R\}.
\end{eqnarray*}
 For any $D_{\lambda}\in Der_{\a^k}(R)$ and $D'_{\mu-\lambda}\in Der_{\a^s}(R)$, define their commutator $[D_\lambda D']_{\mu}$ by
\begin{eqnarray}
[D_\lambda D']_{\mu}(a)=D_\lambda(D'_{\mu-\lambda}a)-(-1)^{|D||D'|}D'_{\mu-\lambda}(D_\lambda a),~~~\forall a\in R.
\end{eqnarray}

\begin{lemma} For any $D_{\lambda}\in Der_{\a^k}(R)$ and $D'_{\mu-\lambda}\in Der_{\a^s}(R)$, we have
\begin{eqnarray*}
[D_\lambda D'] \in Der_{\a^{k+s}}(R)[\lambda].
\end{eqnarray*}
\end{lemma}

{\bf Proof.} For any $a,b\in R$, we have
\begin{eqnarray*}
&&[D_\lambda D']_{\mu}(\partial a)\\
&=&D_{\lambda}(D'_{\mu-\lambda}\partial a)-(-1)^{|D||D'|}D'_{\mu-\lambda}(D_\lambda \partial a)\\
&=& D_\lambda((\partial+\mu-\lambda)D'_{\mu-\lambda}a)+(-1)^{|D||D'|}D'_{\mu-\lambda}((\mu+\lambda)D_\lambda a)\\
&=& (\partial+\mu)D_{\lambda}(D'_{\mu-\lambda}a)-(-1)^{|D||D'|}(\partial+\mu)D'_{\mu-\lambda}(D_{\lambda}a)\\
&=& (\partial+\mu)[D_\lambda D']_{\mu}(a).
\end{eqnarray*}
Also, we have
\begin{eqnarray*}
&&[D_\lambda D']_{\mu}([a_{\g}b])\\
&=& D_\lambda(D'_{\mu-\lambda}[a_{\g}b])-(-1)^{|D||D'|}D'_{\mu-\lambda}(D_\lambda[a_{\g}b])\\
&=&\delta^sD_\lambda([D'_{\mu-\lambda}(a)_{\mu-\lambda+\gamma}\a^s(b)]+(-1)^{|x||D|}[\a^s(a)_{\g}D'_{\mu-\lambda}(b)])\\
&& -\delta^k(-1)^{|D||D'|}D'_{\mu-\lambda}([D_\lambda(a)_{\lambda+\gamma}\a^k(b)]+(-1)^{|x||D|}[\a^k(a)_{\g}D_{\lambda}(b)])\\
&=& \delta^{k+s}[D_\lambda(D'_{\mu-\lambda}(a))_{\mu+\gamma}\a^{k+s}(b)]+\delta^{k+s}(-1)^{|D||D'(x)|}[\a^{k}(D'_{\mu-\lambda}(a))_{\mu-\lambda+\g}D_{\lambda}(\a^{s}(b)))\\
&&+\delta^{k+s}(-1)^{|x||D'|}[D_\lambda(\a^s(a))_{\lambda+\gamma}\a^k(D'_{\mu-\lambda}(b))]+\delta^{k+s}(-1)^{|x||D|}[\a^{k+s}(a)_\g(D_\lambda(D'_{\mu-\lambda}(b)))]\\
&&-\delta^{k+s}(-1)^{|D||D'|}[(D'_{\mu-\lambda}D_\lambda(a))_{\mu+\gamma}\a^{k+s}(b)]+\delta^{k+s}(-1)^{|D'||D(x)|}[\a^{s}(D_\lambda(s))_{\lambda+\gamma}(D'_{\mu-\lambda}(\a^k(b)))]\\
&&-\delta^{k+s}(-1)^{|D|(|D'|+|x|)}[(D'_{\mu-\lambda}(\a^{k}(a)))_{\mu-\lambda+\gamma}\a^{s}(D_{\lambda}(b))]+\delta^{k+s}(-1)^{|x||D'|}[\a^{k+s}(a)_{\lambda}(D'_{\mu-\lambda}(D_\lambda(b)))]\\
&=&\delta^{k+s}[([D_\lambda D']_\mu a)_{\mu+\g}\a^{k+s}(b)]+\delta^{k+s}(-1)^{|x||[D,D']|}[\a^{k+s}(a)_{\g}([D_\lambda D']_{\mu}b)].
\end{eqnarray*}
Therefore, $[D_\lambda D'] \in Der_{\a^{k+s}}(R)[\lambda]$. \hfill $\square$

At the end of this section, we give an application of the $\a$-derivations of a  regular
$\delta$-Hom-Jordan Lie conformal superalgebra $(R, \delta, \a)$. For any $D_\lambda\in Cend(R)$,  define a bilinear
operation $[\c_\lambda \c]_D$ on the vector space $R\oplus \mathbb{R}D$ by
\begin{eqnarray}
[(a+mD)_\lambda(b+nD)]_D=[a_\lambda b]+mD_\lambda(b)-\delta(-1)^{|a||b|} nD_{-\lambda-\partial}(a), \forall a,b\in R, m,n\in \mathbb{R}.~~~~~~
\end{eqnarray}
and a linear map $\a': R\oplus \mathbb{R}D\rightarrow R\oplus \mathbb{R}D$ by $\a'(a+D)=\a(a)+D$.

\begin{proposition}
  $(R\oplus \mathbb{R}D, \a')$ is a  regular $\delta$-Hom-Jordan Lie conformal superalgebra if and only
if $D_\lambda$ is an $\a$-derivation of $(R, \delta, \a)$ with $(1-\delta)D \circ D=0$.
\end{proposition}

{\bf Proof.} Suppose that $(R\oplus \mathbb{R}D, \a')$ is a   regular $\delta$-Hom-Jordan Lie conformal superalgebra. For any $a,b\in R$ and $m,n\in \mathbb{R}$, we have
\begin{eqnarray*}
&&\a'[(a+mD)_\lambda(b+nD)]_D\\
&=&\a'([a_\lambda b]+mD_\lambda(b)-\delta(-1)^{|a||b|}nD_{-\lambda-\partial}(a))\\
&=& \a[a_\lambda b]+m\a (D_\lambda(b))-\delta(-1)^{|a||b|}n\a(D_{-\lambda-\partial}(a)),\\
&&[\a'(a+mD)_{\lambda}\a'(b+nD)]\\
&=&[\a(a)+mD_{\lambda}\a(b)+nD]\\
&=&[\a(a)_\lambda \a(b)]+mD_\lambda \a(b)-\delta (-1)^{|a||b|}nD_{-\lambda-\partial} \a(a).
\end{eqnarray*}
Thus, $\a'$ is an algebra morphism if and only if
\begin{eqnarray*}
\a\circ D_\lambda=D_\lambda\circ \a.
\end{eqnarray*}
Next, if the Hom-Jacobi identity is satisfied, then the following condition holds
\begin{eqnarray*}
D_{\mu}([a_\lambda b])=\delta[(D_\mu a)_{\lambda+\mu}\a(b)]+\delta(-1)^{|a||D|}[\a(a)_{\lambda}(D_\mu b)].
\end{eqnarray*}
Conversely,  if  $D_\lambda$ is an $\a$-derivation of $(R, \delta, \a)$ and $(1-\delta)D \circ D=0$.  For any $a,b\in R$, $m,n\in \mathbb{R}$,
\begin{eqnarray*}
&&[(b+nD)_{-\partial-\lambda}(a+mD)]_D\\
&=& [b_{-\partial-\lambda}a]+nD_{-\partial-\lambda}(a)-\delta(-1)^{|a||b|} mD_\lambda(b)\\
&=&-\delta(-1)^{|a||b|}([a_\lambda b]-(-1)^{|a||b|}\delta nD_{-\partial-\lambda}(a)+mD_\lambda(b))\\
&=& -\delta(-1)^{|a||b|}[(a+mD)_{\lambda}(b+nD)]_D
\end{eqnarray*}
which proves (2.2). And it is obvious that
\begin{eqnarray*}
&&[\partial D_\lambda a]_D=-\lambda[D_\lambda a]_D,\\
&& [\partial a_\lambda D]_D=-D_{-\partial-\lambda}(\partial a)=-\lambda[a_\lambda D]_D,\\
&& [D_\lambda \partial a]_D=D_\lambda(\partial a)=(\partial+\lambda)D_\lambda(a)=(\partial+\lambda)[D_\lambda a]_D,\\
&&[a_\lambda \partial D]_D=-(\partial D)_{-\lambda-\partial}a=(\partial+\lambda)[a_\lambda D]_D,\\
&& \a'\circ \partial =\partial\circ \a'.
\end{eqnarray*}
Thus (2.1) follows. It is easy to  check that  the Hom-Jacobi identity by use of $(1-\delta)D \circ D=0$ and the proof process is left to the reader.   \hfill $\square$

\begin{center}
 {\bf ACKNOWLEDGEMENT}
 \end{center}

    The work of S. J. Guo is supported by  the NSF of China (No. 11761017)
  and the Youth Project for Natural Science Foundation of Guizhou provincial department of education (No. KY[2018]155).
   The work of S. X. Wang is  supported by  the outstanding top-notch talent cultivation project of Anhui Province (No. gxfx2017123)
 and the Anhui Provincial Natural Science Foundation (No. 1808085MA14).

\renewcommand{\refname}{REFERENCES}


\begin{thebibliography}{99}


\bibitem{Ammar2010} F. Ammar,   A. Makhlouf,   Hom-Lie superalgerbas and Hom-Lie asmissible superalgebras,
J. Algebra, 324(2010), 1513-1528.

\bibitem{Ammar2013} F.  Ammar,  A. Makhlouf,  N. Saadaoui,  Cohomology of Hom-Lie superalgebras and $q$-deformed Witt superalgebra,
Czech. Math. J., 63(2013), 721-761.


\bibitem{Kac98} V. G. Kac, Vertex Algebras for Beginners,  Univ. Lect. Ser.  vol. 10 (Amer. Math. Soc., Providence, RI, 1996), second edition 1998.


\bibitem{MC18} L. Ma, L. Chen, J. Zhang, $\delta$-Hom-Jordan Lie superalgebras, Comm. Algebra   46(4)  (2018),  1668-1697.

\bibitem{OK97} S. Okubo, N. Kamiya,   Jordan-Lie superalgebra and Jordan-Lie triple system,  J. Algebra 198(2) (1997),  388-411.


\bibitem{Xu2000}  X. Xu, Quadratic conformal superalgebras, J. Algebra, 231 (2000), 1-38.

\bibitem{Yuan14} L. Yuan,  Hom-Gel'fand-Dorfman bialgebras and Hom-Lie conformal algebras, J. Math. Phys., 55
(2014), 043507 .

\bibitem{Yuan17} L. Yuan, S. Chen,  C. He,  Hom-Gel'fand-Dorfman super-bialgebras and Hom-Lie conformal superalgebras, Acta Mathematica Sinic, English Series, 33(1),
(2017), 96-116.

 \bibitem{Zhao2017} J. Zhao,   L. Chen,  L. Yuan,   Deformations and generalized derivations of
Lie conformal superalgebras,  J. Math. Phys., 58
(2017),  111702.


 \bibitem{Zhao2016} J. Zhao, L. Yuan,  L. Chen,   Deformations and generalized derivations of
Hom-Lie conformal algebras,   Sci. China Math.,  61 (2018), 797-812.

\end{thebibliography}
\end{document}